\documentclass[leqno,12pt]{article}
\usepackage{amssymb,amsfonts}
\usepackage{amsmath,latexsym}

\usepackage{amstext,epic,eepic,epsf,pslatex}
\usepackage{graphicx}

\newcommand{\bepr}{{\em Proof} } 
\newcommand{\enpr}{\hfill \rule{.5em}{.5em}}

\newcommand{\R}{{\mathbb R}}
\newcommand{\N}{{\mathbb N}}

\DeclareMathOperator{\ls}{limsup}

\def\XXint#1#2#3{{\setbox0=\hbox{$#1{#2#3}{\int}$ }
\vcenter{\hbox{$#2#3$ }}\kern-.6\wd0}}

\newtheorem{defin}{Definition}[section] 
\newtheorem{prop}{Proposition}[section] 
\newtheorem{thm}{Theorem}[section] 
\newtheorem{lemma}{Lemma}[section]

\begin{document}

\title{Asymptotic stability of scalar multi-D inviscid shock waves}

\author{Denis Serre \\ \'Ecole Normale Sup\'erieure de Lyon\thanks{U.M.P.A., UMR CNRS--ENSL \# 5669. 46 all\'ee d'Italie, 69364 Lyon cedex 07. France. {\tt denis.serre@ens-lyon.fr}}}

\date{\em In memoriam Andrew J. Majda}

\maketitle

\begin{abstract}
In several space dimensions, scalar shock waves between two constant states $u_\pm$ are not necessarily planar. We describe them in detail. Then we prove their asymptotic stability, assuming that they are uniformly non-characteristic. Our result is conditional for a general flux, while unconditional for the multi-D Burgers equation. 
\end{abstract}



\paragraph{Notations.} The norm in $L^p(\R^d)$ is $\|\cdot\|_p$ and that in $\R^d$ is $|\cdot|$. The positive/negative parts of a real number $s$ are $s^\pm$. The ball of center $x\in \R^d$ and radius $\epsilon$ is denoted $B(x;\epsilon)$. The cone $A(u_-,u_+)$ and its dual $A(u_-,u_+)^\circ$, defined in Section \ref{s:NPshocks}, are used repeatedly in the other parts.

\section{Introduction}

We are concerned with inviscid scalar conservation laws, in arbitrary space dimension $d\ge2$~:
\begin{equation}
\label{eq:scl}
\partial_tu+{\rm div}_xf(u)=0,\qquad t>0,\,x\in\R^d.
\end{equation}
The flux $f:\R\rightarrow\R^d$ is a smooth function. We shall assume the property of non-degeneracy, which we borrow from \cite{LPT}~:
\begin{quote}
{\bf (ND)} For every $(\tau,\xi)\in\R^{1+d}$ such that $\tau^2+|\xi|^2=1$ (or equivalently $\ne0$), 
$${\rm meas}\{s\in\R\,;\,\tau+f'(s)\cdot\xi=0\}=0.
$$
\end{quote}
If $\xi\in\R^d$ is a non-zero vector, the directional flux $s\mapsto\xi\cdot f(s)$ cannot be affine on a non-trivial interval. Remark also that the graph of $f$ may not be locally contained in a hyperplane.
A stronger form of non-degeneracy occurs when the function $s\mapsto \det(f''(s),\ldots,f^{(d+1)}(s))$ does not vanish -- a natural extension of Genuine Nonlinearity to several space dimensions. A paradigm of the latter situation is the flux
$$f^B(s)=(s^2,\ldots,s^{d+1}),$$
which characterizes the so-called {\em multi-D Burgers equation}~:
\begin{equation}\label{eq:Bur}
\partial_tu+\partial_1(u^2)+\cdots+\partial_d(u^{d+1})=0.
\end{equation}

We are interested in the Cauchy problem, where an initial datum
\begin{equation}\label{eq:id}
u(0,x)=a(x),\qquad x\in\R^d,
\end{equation}
is prescribed. Since the seminal work by Kru\v{z}kov \cite{Kru}, it has been known that the notion of solution must be understood in the sense of {\em entropy solutions.} Within this context, the forward Cauchy problem is well-posed for data in $L^\infty(\R^d)$. The existence and uniqueness allow us to define a semi-group $S_t:L^\infty(\R^d)\rightarrow L^\infty(\R^d)$. Its most important properties, which we use below, are the comparison principle and the $L^1$-contraction~: 
\begin{description}
\item[Comp.]
If $a,b\in L^\infty(\R^d)$ are such that $a\le b$ (inequalities are always understood almost everywhere), then for every $t>0$, there holds $S_ta\le S_tb$. In particular, the essential supremum of $S_ta$ is a non-increasing function of time.
\item[Contr.]
If $a,b\in L^\infty(\R^d)$ are such that $b-a\in L^1(\R^d)$, then for every $t>0$, there holds $S_tb-S_ta\in L^1(\R^d)$, and the function $t\mapsto\|S_tb-S_ta\|_1$ is monotonous non-increasing.
\end{description}
The latter property allows us to extend by continuity the semi-group to data $a\in (L^\infty+L^1)(\R^d)$. Whether the corresponding flow $u(t):=S_ta$ is a solution in Kru\v{z}kov's sense is still an open question~; this is discussed in \cite{DSLS} in the case of (\ref{eq:Bur}).
The reader who is familiar with Kru\v{z}kov's  theory knows that the properties above are associated with the conservation of the mass (whence the terminology of conservation law),
$$\int_{\R^d}(S_tb-S_ta)\,dx\equiv\int_{\R^d}(b-a)\,dx.$$
This information will however be of little help in the multi-dimensional stability analysis. In one space dimension instead, it is used to determine the shift between the shocks at initial and final times, in terms of the mass of the initial disturbance.

\bigskip

The simplest solutions of (\ref{eq:scl}), besides the constants, are planar waves $u(t,x)=U(x\cdot\nu-\sigma t)$, where $\nu\in S_{d-1}$ is the direction of propagation, and $\sigma\in\R$ is the normal velocity. The profile $U$ obeys to $f(U)\cdot\nu-\sigma U={\rm cst}$, which implies, because of the non-degeneracy and the entropy condition, that $U$ is piecewise constant. Generically, $U$ takes only two values:
$$U(y)=\left\{\begin{array}{lrc}
u_-, & \hbox{if} & y<y^*, \\
u_+, & \hbox{if} & y>y^*,
\end{array}\right.$$
and we speak of a {\em planar shock wave}. Up to the flip $(u_-,u_+;\nu,\sigma)\leftrightarrow(u_+,u_-;-\nu,-\sigma)$, we may always assume that $u_+<u_-$. Then the admissibility criterion given by the entropy condition writes
\begin{equation}\label{eq:Ol}
f(s)\cdot\nu-\sigma s\le f(u_\pm)\cdot\nu-\sigma u_\pm,\qquad\forall s\in(u_+,u_-).
\end{equation}
Recall that by letting $s\rightarrow u_\pm$, (\ref{eq:Ol}) implies the Lax shock condition
\begin{equation}\label{eq:Laxplan}
f'(u_+)\cdot\nu\le\sigma \le f'(u_-)\cdot\nu.
\end{equation}

\bigskip

A natural question is whether planar shock waves are stable under localized disturbances: 
\begin{quote}
{\em Given an initial data $a(x)=U(x\cdot\nu)+\phi(x)$, where $\phi\in{\cal D}(\R^d)$, how does the solution $u(t)=S_ta$  behave as $t\rightarrow+\infty$~?} We shall relax the restrictions on both $U$ and $\phi$ in the statements below.
\end{quote}
A naive, though deadly false, guess would be that $u(t)-U(x\cdot\nu-\sigma t)\rightarrow0$ in some sense, perhaps up to a shift (orbital stability). We shall see that, whenever $d\ge2$, the answer is not so simple. As described in Section \ref{s:NPshocks}, a planar shock is just  one among all the shocks between $u_-$ and $u_+$, without any special property. Their set can be split into equivalence classes in which two shock profiles differ only by an integrable disturbance, or by a compactly supported one. Our main result is that such classes are stable, a fact that can be interpreted as a generalized form of orbital stability. Mind that since $d\ge2$, a shock profile and its shifts belong  in general to distinct equivalence classes.

Since we shall exploit the $L^1$-contraction property of the Kru\v{z}kov semigroup, we understand the stability  in terms of an $L^1$-distance. Our most general result thus concerns disturbances $\phi$ that belong to $L^1(\R^d)$. Notice that because this perturbation can be arbitrarily wild, the shock structure is lost in the transient times, even though we recover it in the time asymptotics. Thus we shall not try to describe the evolution of the perturbed shock front. 

\bigskip

We start our analysis by a description of general shock waves between constant states, for which we refer to Theorem \ref{th:struct} below. 
An important fact is that the shock fronts are Lipschitz graphs, without higher regularity in general. This lack of regularity makes it difficult, if not impossible, to prove their local-in-time stability in the spirit of A. Majda's memoirs \cite{Maj1,Maj2}. This flaw is closely related to the fact that, because $d\ge2$, the Lopatinski\u{\i} condition for a scalar shock is satisfied only in a non-uniform way. 

The other results address the asymptotic stability of such shocks, when they are uniformly non-characteristic. The first one deals with the case where the data take values in the interval $[u_+,u_-]$ defined by the end states~; it is valid for arbitrary fluxes. For the sake of simplicity, the shocks under consideration are stationary, in the sense that $f(u_+)=f(u_-)$. This does not limit the generality, since we may always assume this situation by choosing an appropriate moving frame, which amounts to adding a linear function to the flux.
\begin{thm}\label{th:upm}
We assume ({\bf ND}). Let $U=U(x)$ be a steady shock between two constant states $u_+<u_-$ (for a general description, see Theorem \ref{th:struct}). We assume that the shock front is uniformly non-characteristic, in the sense of Definition \ref{d:UNC}. Let $\phi(x)\in L^1(\R^d)$ be an initial perturbation, such that the corresponding initial data $a=U+\phi$ still takes values in $[u_+,u_-]$ (in particular, $\phi$ is bounded).

Then the solution $u(t)=S_ta$ tends in the $L^1$-distance towards another steady shock $\hat U$ in the same class as $U$ {\em modulo} $L^1(\R^d)$~:
\begin{eqnarray*}
\hat U :\,\R^d & \mapsto & \{u_-,u_+\}, \\
\hat U-U & \in & L^1(\R^d), \\
\|u(t)-\hat U\|_1 & \stackrel{t\rightarrow+\infty}\longrightarrow & 0.
\end{eqnarray*}
\end{thm}

The proof given in Section \ref{s:upm},  reminiscent to dynamical systems theory, follows the strategy developed in \cite{Ser_HB}. We start with the subcase where the initial perturbation is compactly supported, and we prove that it remains so, uniformly in time. The kinetic theory of Lions, Perthame \& Tadmor \cite{LPT} ensures the compactness of the sequence $(u(\cdot+\tau,\cdot))_{\tau\rightarrow+\infty}$ and allows us to consider the $\omega$-limit set $\Omega$. The contraction property provides us with a large family of Lyapunov functions $u\mapsto\|W-u\|_1$, indexed by the shock waves $W$ that belong to the same class as $U$. An element $\bar u\in\Omega$ is a solution of (\ref{eq:scl}), still taking values in $[u_+,u_-]$, with the special property that each of these Lyapunov functions remains constant as the time evolves (Lasalle's Invariance Principle). We show that this implies that $\bar u$ is itself a steady shock in the same class.

Eventually, we remove the assumption that $\phi$ be compactly supported, by using a standard density\,/\,contraction argument.

\bigskip

When the initial data is not confined to the interval $[u_+,u_-]$, a reasonnable strategy is to prove that the ``overhead'' $\{x\,;\,u(t,x)>u_-\}$ disappears in finite time, as well as its counterpart $\{x\,;\,u(t,x)<u_+\}$. If so, we may apply Theorem \ref{th:upm} above to the data $u(T)$ for some $T>0$ large enough. We achieve this goal in the case of the multi-D Burgers equation:
\begin{thm}\label{th:over}
Consider the multi-D Burgers equation ($f=f^B$) and a (non-necessarily planar) shock $U(x-tv)$ between the states $u_+<u_-$, whose shock front is uniformly non-characteristic. Let $\phi(x)\in L^\infty(\R^d)$ be a compactly supported initial perturbation. Let us form the initial data $a(x)=U(x)+\phi(x)$.

Then there exists a finite time $T$ such that, for every $t>T$, the solution $u(t)=S_ta$ takes values in $[u_+,u_-]$.
\end{thm}

A uniform decay, at an algebraic rate, of the overhead is established as a consequence of a dispersion property, established by L. Silvestre and the author \cite{DSLS}, given in the Appendix. Once its amplitude is small enough, the non-characteristicness forces it to move toward the shock, where it is absorbed.

Notice that if the shock was characteristic, for instance if $f'(u_-)\cdot\nu=0$ in a direction $\nu$ normal to the shock, while $f''(u_-)\cdot\nu\ne0$, then the statement would fail. A small amplitude overhead would travel away from the shock and its mass could never be absorbed. This phenomenon occurs already in one-dimensional case, as described in Paragraph 3.4 of \cite{Ser_HB}. It is specific to the case of inviscid conservation law and does not happen when some amount of dissipation is added to the equation (see \cite{FrS,Ser_HB} for $d=1$ and \cite{GM,HZ} for $d=2$).

\bigskip

Combining the statements above, we obtain
\begin{thm}\label{th:Bur}
Consider the multi-D Burgers equation ($f=f^B$) and a uniformly non-characteristic shock $U(x-tv)$ between the states $u_+<u_-$. Let $\phi(x)\in  L^1(\R^d)$ be an initial perturbation, with which we form the initial data $a(x)=U(x)+\phi(x)$.
\begin{enumerate}
\item Then the solution $u(t)=S_ta$ tends in the $L^1$-distance towards another shock in the same class as $U$ {\em modulo} $L^1(\R^d)$~:
\begin{eqnarray*}
\hat U :\,\R^d & \mapsto & \{u_-,u_+\}, \\
\hat U-U & \in & L^1(\R^d), \\
\|u(t)-\hat U(\cdot-tv)\|_1 & \stackrel{t\rightarrow+\infty}\longrightarrow & 0.
\end{eqnarray*}
\item If $\phi\in L^\infty(\R^d)$ is compactly supported, then so is $\hat U-U$.
\end{enumerate}
\end{thm}

We emphasize the generality of our stability results. Neither the jump $u_+-u_-$, nor the size of the amplitude or the mass of the perturbation $\phi$ need to be small. The width between a necessary condition -- the non-characteristicness -- and our sufficient one is very slim. The counterpart is that we do not know at which rate the convergence takes place. In this respect, our results differ significantly from those of Duch\^ene \& Rodrigues \cite{DR1}, where the stability is induced by a lower-order term (balance law) and the rate is exponential. We cannot either characterize the limit $\hat U$ in terms of the initial data~; we can only say that
$$\int_{\R^d}(\hat U(x)-U(x))\,dx=\int_{\R^d}\phi(x)\,dx.$$

\bigskip

Remark that the non-characteristicness (in a uniform manner) is needed only away from a compact subset. As a matter of fact, we may always modify the profile $U$ in a bounded set, where we can easily let it uniformly non-characteristic. 

\paragraph{Plan of the paper.} Shocks taking values in a doubleton $\{u_-,u_+\}$ are described in Section \ref{s:NPshocks}. Special attention is given to non-characteristic shocks. We prove Theorem \ref{th:upm} in Section \ref{s:upm}, and Theorem \ref{th:over} in Section \ref{s:over}. Two anterior results of the theory of scalar conservation laws are listed in Appendix \ref{ap}. 

\paragraph{Remark.} Almost everything in this paper is specific to the scalar case. Not only we use extensively the comparison principle and the $L^1$-contraction, the latter being a source of Lyapunov functions. But on another hand, a pair $(u_-,u_+)$ of constant states yields a shock wave in every direction taken in some (usually wide) cone, denoted below $A(u_-,u_+)^\circ$. This contrasts with the system case, where not all pairs can be related by a single shock~; and even when they can be, it is usually in only one direction. For instance, two states of full gas dynamics with density $\rho_\pm$, pressure $p_\pm$ and specific internal energy $e_\pm$ must satisfy the compatibility condition
$$e_+-e_-+\frac{p_++p_-}2\,\left(\frac1{\rho_+}-\frac1{\rho_-}\right)=0$$
in order that a shock between them be possible~; and then the normal to the shock must be colinear to the jump $\vec u_+-\vec u_-$ of the fluid velocity.

\section{Non-planar shock waves}\label{s:NPshocks}

Recall that $f=(f_1,\ldots,f_d):\R\rightarrow\R^d$ is a smooth flux, which satisfies the non-degeneracy assumption ({\bf ND}).

Given two real numbers $u_+<u_-$, we define the velocity
$$v(u_-,u_+):=\frac1{u_+-u_-}\,(f(u_+)-f(u_-)).$$
If $\xi\in S_{d-1}$, the Rankine-Hugoniot condition in the direction $\xi$ provides us with a number
$$\sigma(u_-,u_+,\xi):=\xi\cdot v(u_-,u_+),$$
which is the normal velocity of a discontinuity $u_-\mapsto u_+$ along a hyperplane $\xi\cdot x={\rm cst}$. Then we extend the definition of $\sigma$ to every $\xi\in\R^d$. 

We say that $\xi$ is an admissible direction if $\xi\cdot f$ satisfies also the Oleinik condition
\begin{equation}\label{eq:Olxi}
\xi\cdot f(s)-\sigma(u_-,u_+;\xi)s\le\xi\cdot f(u_\pm)-\sigma(u_-,u_+;\xi)u_\pm,\qquad\forall s\in(u_+,u_-).
\end{equation}
This is equivalent to saying that (\ref{eq:scl}) admits the entropy solution $U(x\cdot\xi-\sigma t)$, where $U(y)\equiv u_\pm$, according to whether $y$ is positive or negative. If $\xi\ne0$, this is a planar shock wave in the direction $\xi$. Of course, $\xi=0$ is admissible. 

\bigskip

The set of the admissible directions, denoted $A(u_-,u_+)$, is a closed convex cone. If it reduces to $\{0\}$, then a shock $u_-\mapsto u_+$ is not possible whatever the direction. Because of the non-degeneracy, $A(u_-,u_+)$ cannot contain an entire line.

\bigskip

Suppose now that $A(u_-,u_+)$ is non-trivial. We are interested in all the solutions $u$ of (\ref{eq:scl}) which take only the values $u_\pm$. Let us define the auxiliary flux $F(s):=f(s)-sv(u_-,u_+)$, for which we have $F(u_+)=F(u_-)$. This common value is denoted $\bar F$. For the sake of simplicity, we rewrite the conservation law in  the moving frame $(t,x':=x-tv(u_-,u_+))$, where it becomes
\begin{equation}\label{eq:mf}
\partial_tu+{\rm div}\,F(u)=0.
\end{equation}
Our assumption that $u(t,x)\equiv u_\pm$ a.e. implies $F(u)\equiv\bar F$. The conservation law thus reduces to $\partial_tu=0$~: the solution is steady in the moving frame. From now on, we write $x$ instead of $x'$ and $u(t,x)=u(x)$. Let us point out that if we replace $f$ by $F$ in the definition above, then $\sigma\equiv0$.

The equation being understood in Kru\v{z}kov's sense, we must also write the entropy inequalities, namely
$$\partial_t|u-k|+{\rm div}\,({\rm sgn}(u-k)(F(u)-f(k)))\le0,$$
for every real parameter $k$. Since $u$ is stationary and $F(u)$ is constant, this reduces to
$$(\bar F-F(k))\cdot\nabla{\rm sgn}(u-k)\le0.$$
When either $k>u_-$ or $k<u_+$, this inequality is trivial since $u-k$ is of constant sign. For $k\in(u_+,u_-)$ instead, it writes $(\bar F-F(k))\cdot\nabla\chi\ge0$ where $\chi(x)=\pm1$ according to whether $u(x)=u_\pm$. This amounts to saying that $\nabla\chi$ takes values in $A(u_-,u_+)$. The function $\chi$ is thus non-decreasing in the directions of the dual cone
$$A(u_-,u_+)^\circ=\left\{\vec n\in\R^d\,;\,\vec n\cdot\xi\ge0,\forall \xi\in A(u_-,u_+)\right\}.$$
Remark that $A(u_-,u_+)^\circ$ is nothing but the convex cone spanned by the vectors $\bar F-F(s)$ as $s\in[u_+,u_-]$.
We infer that
\begin{lemma}\label{l:cones} 
If $u$ equals $u_+$ at some Lebesgue point $x$, then $u\equiv u_+$ in the cone $x+A(u_-,u_+)^\circ$. Likewise, if $u(x)=u_-$, then $u\equiv u_-$ in $x-A(u_-,u_+)^\circ$.  
\end{lemma}

Let us pick a vector $W$ in the interior of $A(u_-,u_+)^\circ$, which exists because of ({\bf ND}). Suppose that $u$ is not a constant function. Thus there exists two points $x_\pm$ such that $u(x_\pm)=u_\pm$. Let $L$ be any line of direction $W$. Then $L\cap(x_++A(u_-,u_+)^\circ)\ne\emptyset$ and therefore $L$ contains a point at which $u=u_+$. Symmetrically, $L\cap(x_--A(u_-,u_+)^\circ)\ne\emptyset$ and $L$ must contain a point at which $u=u_-$. Thus there exists a point $X(L)\in L$ which divides $L$ into two half-lines on which $u\equiv u_\pm$, respectively. We now pick an element $\lambda$  of $A(u_-,u_+)$. The cone $A(u_-,u_+)^\circ$ is contained in a half-space delimited by the hyperplane $H$ of equation $\lambda\cdot x=0$. This allows us to define the coordinates $(y,r)\in H\times\R$ associated with the decomposition $\R^d=H\oplus\R W$. The set of lines $L$ as above is parametrized by $y$, and the divider writes $X(L)=y+\psi(y)W$ for some function $\psi$. The shock front is therefore the {\em graph} of the function $\psi$ above $H$. And because of Lemma \ref{l:cones}, we have $|\psi(y')-\psi(y)|\le\psi_0(y'-y)$, where $\psi_0$ is the function whose epigraph is $A(u_-,u_+)^\circ$. Since $W$ is interior to $A(u_-,u_+)^\circ$, $\psi_0$ satisfies an inequality $\psi_0(y)\le C|y|$ and thus $\psi$ is globally Lipschitz.

We summarize our analysis in the following statement.
\begin{thm}\label{th:struct}
We assume ({\bf ND}).
Suppose that $A(u_-,u_+)\ne\{0\}$. Then every entropy solution $u$ of (\ref{eq:scl}) which takes values in $\{u_-,u_+\}$ is stationary and is of one of the following forms:
\begin{description}
\item[Constants.] Either $u\equiv u_-$ in $\R^d$, or $u\equiv u_+$ in $\R^d$.
\item[Shocks.] A complete Lipschitz hypersurface $\Gamma$, whose normals belong to $A(u_-,u_+)$, separates $\R^d$ into two halves $D_\pm$, on which $u\equiv u_\pm$ respectively. In particular $\Gamma$ is a graph in suitable linear coordinates.
\end{description}
In the latter case, for every interior point $W$ of $A(u_-,u_+)^\circ$, the family $(D_-+sW)_{s\in\R}$ is ordered by inclusion and we have
\begin{equation}\label{eq:capcup}
\bigcap_s(D_-+sW)=\emptyset,\qquad\bigcup_s(D_-+sW)=\R^d.
\end{equation}
\end{thm}
Notice that in the limit case where $A(u_-,u_+)=\R_+\lambda$ for some $\lambda\in S_{d-1}$, then the non-constant solutions described in the theorem above are rigid: they are planar shock waves $U(\lambda\cdot x-t\sigma(u_-,u_+;\lambda))$.

\subsection{Characteristic shocks}

Recall that Lipschitz hypersurfaces admit a normal direction at almost every point. The following definition deviates slightly from standards. We state it for globally defined shock profiles, though it can be used also for local ones.
\begin{defin}\label{d:NC}
Consider a shock as described in Theorem \ref{th:struct}. Let $x\mapsto\nu(x)$ be the measurable map defined over the shock front, where $\nu(x)$ is the unit normal, oriented towards the domain $D_+$. Recall that $\nu(x)\in A(u_-,u_+)$. 

We say that the shock is {\em characteristic} at $x$ if $\nu(x)$ belongs to the boundary of $A(u_-,u_+)$. It is {\em non-characteristic} otherwise, that is when $\nu(x)$ is an interior point of the cone.
\end{defin}

Two distinct phenomena may cause an admissible shock with normal direction $\nu\in S_{d-1}$ to be characteristic:
\begin{itemize}
\item One of the Lax inequalities (\ref{eq:Laxplan})  is non-strict: either $\nu\cdot F'(u_+)=0$, or $\nu\cdot F'(u_-)=0$.
\item The (non-positive) function $s\mapsto\nu\cdot (F(s)-\bar F)$ vanishes at some interior point $\bar s\in(u_+,u_-)$.
\end{itemize}
The former possibility is the one that most authors consider usually. It forbids the $L^1$-asymptotic stability of the corresponding shock, because it allows some charateristic lines to emerge tangentially from the shock. Some non-trivial mass in excess (above $u_-$ or below $u_+$) can escape as the time increases. We think that the second possibility is equally important, because then a general initial disturbance lets the shock split into two subshocks $u_-\mapsto \bar s$ and $\bar s\mapsto u_+$, each one being characteristic in the ordinary sense. 

We wish here to avoid both difficulties, whence the definition above. But in practice, because the shock front is not necessarily smooth, and also that it extends to infinity, we need some uniformity:
\begin{defin}\label{d:UNC}
With the same notations as in Definition \ref{d:NC}, we say that the shock is {\em uniformly} non-characteristic if $A(u_-,u_+)$ has a non-empty interior, and if the Gau\ss\, map $x\mapsto\nu(x)$ takes its values in a compact subset of this interior.
\end{defin}

Obviously, a planar shock wave is uniformly non-characteristic if and only if it is non-characteristic.

The shock front $\{x\,;\,x=y+\psi(y)W\}$ is non-characteristic at a point $x_0=y_0+\psi(y_0)W$ when there exists a constant $0<\rho<1$ such that the function $\psi$ satisfies the enhanced Lipschitz condition
$$|\psi(y)-\psi(y_0)|\le\rho\,\psi_0(y-y_0)$$
in a neighbourhood of $y_0$. It is uniformly non-characteristic when one may choose the same constant $0<\rho<1$ almost everywhere. Since $\psi_0$ is subadditive (convex and homogeneous of degree $1$), uniform non-characteristicness can be rewritten equivalently
\begin{equation}\label{eq:LipPsi}
\exists\rho<1\quad\hbox{s.t.}\quad|\psi(y')-\psi(y)|\le\rho\,\psi_0(y'-y),\qquad\forall\,y,y'\in H.
\end{equation}
The fact that the cone $A(u_-,u_+)$ depends continuously upon its arguments implies that the set of uniformly non-characteristic shock profiles with a given front is open:
\begin{prop}\label{p:autrechoc}
Let the shock $U(x-tv(u_-,u_+))$, between the constant states $u_\pm$,\ be uniformly non-characteristic. Then there exists an $\eta>0$ such that, if $\hat u_\pm$ are chosen according to $|\hat u_\pm-u_\pm|<\eta$, then the function $\hat U$ defined by
$$\hat U(x)=\left\{\begin{array}{lcr}
\hat u_- & \hbox{if} & x\in D_-, \\
\hat u_+ & \hbox{if} & x\in D_+
\end{array}\right.$$
is still the profile of a uniformly non-characteristic shock $\hat U(x-t\hat v)$.
\end{prop}
We point out that $\hat U$ and $U$ have the same shock front though different end states. Of course, the modified velocity $\hat v$, given by the Rankine--Hugoniot condition, is
\begin{equation}\label{eq:hatv}
\hat v=v(\hat u_-,\hat u_+)=\frac1{\hat u_+-\hat u_-}\,(f(\hat u_+)-f(\hat u_-))=v(u_-,u_+)+O(\eta).
\end{equation}

\bigskip

\begin{prop}\label{p:trans}
We assume that the interior of $A(u_-,u_+)$ is not empty.
Let $u$ be a shock from $u_-$ to $u_+$, as described in Theorem \ref{th:struct}. Assume that it is uniformly non-characteristic. Denote $D_\pm$ the domains $\{x\in\R^d\,;\,u(x)=u_\pm\}$. Let $x\in\R^d$ be given.

Then the domains
$$D_-\bigcap(x+A(u_-,u_+)^\circ),\qquad D_+\bigcap(x-A(u_-,u_+)^\circ)$$
are bounded.
\end{prop}

\bepr

It is enough to prove that the first intersection is bounded. In terms of the coordinates $(r,y)\in\R\times H$, $D_-$ is given by $r<\psi(y)$. In particular it is contained in the domain defined by $r<\psi(0)+\rho\psi_0(y)$. On the other hand the cone $x+A(u_-,u_+)^\circ$ has equation $r-r_0\ge\psi_0(y-y_0)$, where $x=:r_0W+y_0$. A point in the intersection thus satisfies
\begin{eqnarray*}
\psi_0(y) & \le & \frac1{1-\rho}\,(\psi(0)+\psi_0(y_0)-r_0), \\
r_0-\psi_0(y_0)\le & r & \le\frac1{1-\rho}\,(\psi(0)+\rho(\psi_0(y_0)-r_0)).
\end{eqnarray*}
Because the interior of $A(u_-,u_+)$ is not empty, $\psi_0$ is $>0$ away from the origin, and the first line above tells us that $y$ belongs to a bounded set of $H$. The second line controls $r$.

\enpr

\section{Proof of Theorem \protect\ref{th:upm}}\label{s:upm}


Let us assume in a first instance that $\phi$ is compactly supported. Because the interior of $A(u_-,u_+)^\circ$ is not empty, there exist two points $x^\pm$ such that 
$${\rm Supp}\,\phi\subset (x_-+A(u_-,u_+)^\circ)\bigcap (x_+-A(u_-,u_+)^\circ).$$
Let us define a function 
$$w(x):=\left\{\begin{array}{lcl}
u_-, & \hbox{if} & x\in D_-\cup(x_+-A(u_-,u_+)^\circ), \\
u_+, & \hbox{if} & \hbox{not.} \end{array}\right.$$
Likewise we define
$$z(x):=\left\{\begin{array}{lcl}
u_+, & \hbox{if} & x\in D_+\cup(x_-+A(u_-,u_+)^\circ), \\
u_-, & \hbox{if} & \hbox{not.} \end{array}\right.$$
Then $z\le a\le w$. Since both $w$ and $z$ are the profiles of steady solutions of (\ref{eq:scl}), we infer (comparison principle) that our solution satisfies
\begin{equation}\label{eq:cadre}
w(x)\le u(t,x)\le z(x),\qquad\forall t>0,\,x\in\R^d.
\end{equation}
In other words, the perturbation $u(t)-U$ remains compactly supported, uniformly in time:
$${\rm Supp}(u(t)-U)\subset (x_-+A(u_-,u_+)^\circ)\bigcap (x_+-A(u_-,u_+)^\circ),\qquad\forall t>0.$$

\bigskip

We now define the time-shifts of the solution, indexed by $\tau>0$~:
$$u^\tau(t,x):=u(t+\tau,x),\qquad (t,x)\in(-\tau,+\infty)\times\R^d.$$
Each $u^\tau$ is an entropy solution of (\ref{eq:scl}). The sequence $(u^\tau)_{\tau>0}$ being uniformly bounded, is is relatively compact in $L^1_{\rm loc}((A,+\infty)\times\R^d)$ for every $A\in\R$ according to Theorem 3\footnote{The cited result makes in addition the spurious assumption that the sequence be bounded in $L^\infty_t(L^1(\R^d))$. Truncating $u^\tau(T)$ to some ball $B$ and using the finite velocity of waves, we infer the compactness within a cone of basis $\{T\}\times B$ and fixed slopes. Such cones can be taken arbitrarily large.} of \cite{LPT}, thanks to Assumption ({\bf ND}). If we now consider the sequence $u^\tau-U$, which is compactly supported in space, the compactness holds instead in $L^1((A_1,A_2)\times\R^d)$ for every finite $A_1<A_2$.

Our next target is to characterize the $\omega$-limit set $\Omega$ of this sequence as $\tau\rightarrow+\infty$. It is made of functions $\bar u$ that are limits of some subsequences $u^{\tau_k}$ (with $\tau_k\rightarrow+\infty$) in  $L^1_{\rm loc}(\R^{1+d})$. Because the convergence is strong, one can pass to the limit in the equation and in the entropy inequalities. We thus find that such a limit is itself an entropy solution of (\ref{eq:scl}), though in the entire space $\R^{1+d}$.

We now use the fact that (\ref{eq:scl}) admits plenty of Lyapunov functions. Suppose that $R(x)$ is a shock from $u_-$ to $u_+$, which coincides with $U$ away from a bounded set. In other words, the shock front coincides with that of $U$ away from a bounded set. Then $a-R\in L^1(\R^d)$, thus $u(t)-R\in L^1(\R^d)$ for all times, and $t\rightarrow\|u(t)-R\|_1$ is a non-increasing function.. Let $\ell_R$ be its limit as $t\rightarrow+\infty$. Then 
\begin{equation}\label{eq:ellW}
\sup_{t\in[A_1,A_2]}\|u^{\tau_k}(t)-R\|_1=\|u(A_1+\tau_k)-R\|_1\rightarrow\ell_R.
\end{equation}
We infer that any $\bar u\in\Omega$ satisfies
\begin{equation}\label{eq:Lyap}
\|\bar u(t)-R\|_1\equiv\ell_R,\qquad t\in\R.
\end{equation}

Recall that the decay of $\|\bar u(t)-R\|_1$ follows from integrating in space the inequality (that Kru\v{z}kov proved as a consequence of entropy inequalities) between two solutions
\begin{equation}\label{eq:entuW}
\partial_t|\bar u-R|+{\rm div}\,[{\rm sgn}(\bar u-R)\,(f(\bar u)-f(R))]\le0.
\end{equation}
Thus the constancy (\ref{eq:Lyap}) tells us that actually (\ref{eq:entuW}) is an equality:
\begin{equation}\label{eq:equW}
\partial_t|\bar u-R|+{\rm div}\,[{\rm sgn}(\bar u-R)\,(f(\bar u)-f(R))]=0.
\end{equation}

We now exploit the flexibility given by the very large set of shocks $R$, to prove the following.
\begin{lemma}\label{l:flex}
For every test function $\theta=\theta(x)\in C_K(\R^d)$, one has
\begin{equation}\label{eq:PDEth}
\partial_t\left[\theta\left(\bar u-\frac{u_-+u_+}2\right)\right]+{\rm div}\,[\theta\,(f(\bar u)-\bar f)]=0,
\end{equation}
where $\bar f$ denotes the common value $f(u_\pm)$.
\end{lemma}

\bepr

Let $\rho\in{\cal D}(\R^{1+d})$ be a test function. The PDE (\ref{eq:equW}) tells us that
$$\int_{\R^{1+d}}\left(|\bar u-R|\partial_t\rho+[{\rm sgn}(\bar u-R)\,(f(\bar u)-f(R))]\cdot\nabla_x\rho\right)\,dx\,dt=0.$$
Because $R$ takes only the values $u_\pm$, this rewrites as
\begin{equation}\label{eq:onzetest}
\int_{\R^{1+d}}\left(|\bar u-R|\partial_t\rho+[{\rm sgn}(\bar u-R)\,(f(\bar u)-\bar f)]\cdot\nabla_x\rho\right)\,dx\,dt=0.
\end{equation}
Let $R'$ be another multi-D shock between the constants $u_\pm$ coinciding with $U$ away from some bounded domain. Thus $R$ and $R'$ coincide away from a bounded domain $D$. Up to replacing $R,R'$ by $\min(R,R')$ and $\max(R,R')$ respectively, we may assume $R\le R'$, which means that $R'\equiv u_-$ while $R\equiv u_+$ within $D$ (recall that $u_+\le u\le u_-$). Writing (\ref{eq:onzetest}) for $R$ and $R'$ respectively, and taking the difference between both equalities, we obtain
\begin{equation}\label{eq:surD}
\int_D\left((2\bar u-u_--u_+)\partial_t\rho+2(f(\bar u)-\bar f))\cdot\nabla_x\rho\right)\,dx\,dt=0.
\end{equation}

The identity (\ref{eq:surD}) is rather flexible: given two points $y_\pm\in\R^d$, one may take
\begin{equation}\label{eq:dia}
D=(y_-+A(u_-,u_+)^\circ)\bigcap(y_+-A(u_-,u_+)^\circ
\end{equation}
by the following construction:
\begin{eqnarray*}
\{x\,;\,R(x)=u_-\} & = & (D_-\bigcup(y_+-A(u_-,u_+)^\circ)\setminus(y_-+A(u_-,u_+)^\circ), \\
\{x\,;\,R'(x)=u_-\} & = & (D_-\setminus(y_-+A(u_-,u_+)^\circ)\bigcup(y_+-A(u_-,u_+)^\circ).
\end{eqnarray*}
We have used above the set theoretic property
$$(B\cup S)\setminus T\subset(B\setminus T)\cup S,\qquad\left((B\setminus T)\cup S\right)\setminus\left((B\cup S)\setminus T\right)=S\cap T,$$
together with the fact that both fronts are still Lipschitz hypersurfaces with normals in $A(u_-,u_+)$, so that both $R$ and $R'$ are admissible shocks. Thanks to Proposition \ref{p:trans}, both $R$ and $R'$ are compact perturbations of $U$.
\begin{quote}
Remark. {\em In terms of the fronts, which are graphs of functions $b,h,k$ respectively, we are forming the functions
$$\hat h(y)=\min\{h(y),\max\{b(y),k(y)\}\},\qquad \hat k(y)=\max\{k(y),\min\{b(y),h(y)\}\},$$
and we have $\hat k-\hat h=(k-h)^+$.}
\end{quote}

Chosing above $y_+$ in the interior of $A(u_-,u_+)^\circ$, and $y_-:=-y_+$, we obtain a specific domain $\Delta$ which is a balanced convex neighbourhood of the origin. It is the unit ball for a suitable norm of $\R^d$. Then appropriate choices of $y_\pm$ yield domains $D=D(y_0,\mu)=y_0+\mu\Delta$ (where $y_0\in\R^d$ and $\mu>0$) that are arbitrary balls for this norm. Writing (\ref{eq:surD}) in the form
$$\int_D m(x)\,dx=0,\qquad m(x):=\int_\R\left((2\bar u-u_--u_+)\partial_t\rho+2(f(\bar u)-\bar f))\cdot\nabla_x\rho\right)\,dt,$$
we obtain that $m$ vanishes at every Lebesgue point, hence almost everywhere. Multiplying now by $\theta$ and integrating over $\R^d$, we obtain
$$\int_{\R^{1+d}}\theta(x)\left((2\bar u-u_--u_+)\partial_t\rho+2(f(\bar u)-\bar f)]\cdot\nabla_x\rho\right)\,dx\,dt=0,$$
which is the weak formulation of (\ref{eq:PDEth}).

\enpr

\bigskip

To conclude, we develop (\ref{eq:PDEth}) into
$$\theta(\partial_t\bar u+{\rm div}\,f(\bar u))+(f(\bar u)-\bar f)\cdot\nabla\theta=0.$$
Because $\bar u$ solves the conservation law, there remains
$$(f(\bar u)-\bar f)\cdot\nabla\theta=0,\qquad\forall\,\theta\in C_K(\R^d).$$
This is equivalent to saying that $f(\bar u)\equiv\bar f$. Remembering that $\bar u(t,x)\in[u_+,u_-]$ and that $f(s)\ne \bar f=f(u_\pm)$ for $s\in(u_+,u_-)$, we infer $\bar u(t,x)\in\{u_-,u_+\}$ almost everywhere. According to Theorem \ref{th:struct}, $\bar u$ is a stationary shock $\hat U(x)$.

So let $\tau_k\rightarrow+\infty$ be such that $u^{\tau_k}$ converges towards $\hat U$. Applying (\ref{eq:ellW}) to $\hat U$, we find  $\ell_{\hat U}=0$, that is
$$\lim_{t\rightarrow+\infty}\|u(t)-\hat U\|_1=0.$$
This proves the convergence of $u(t)$ towards the steady shock $\hat U$.

\bigskip

In the general case, where the perturbation is integrable but not compactly supported, we approach $\phi$ in $L^1(\R^d)$ by a sequence $(\phi_m)_{m\in\N}$ of compactly supported perturbations, such that $a_m:=U+\phi_m$ still take values in $[u_+,u_-]$. Applying the Theorem, already proved under this restriction, to the data $a_m$, we have
$$\lim_{t\rightarrow+\infty}\|S_ta_m-Z_m\|_1=0,$$
for some steady shock $Z_m$ such that $Z_m-U\in L^1(\R^d)$. By the contraction property, we have
$$\|Z_m-Z_p\|_1=\lim_{t\rightarrow+\infty}\|S_ta_m-S_ta_p\|_1\le\|a_m-a_p\|_1=\|\phi_m-\phi_p\|_1\stackrel{m,p\rightarrow+\infty}\longrightarrow0.$$
The sequence $(Z_m)_{m\in\N}$ is thus Cauchy, hence convergent in $U+L^1(\R^d)$. Let $\hat U$ be its limit. From
\begin{eqnarray*}
\|u(t)-\hat U\|_1 & \le & \|u(t)-S_ta_m\|_1+\|S_ta_m-Z_m\|_1+\|Z_m-\hat U\|_1 \\
& \le & 2\|\phi-\phi_m\|_1+\|S_ta_m-Z_m\|_1,
\end{eqnarray*}
we infer
$$\ls_{t\rightarrow+\infty}\|u(t)-\hat U\|_1\le2\|\phi-\phi_m\|_1.$$
Passing to the limit as $m\rightarrow+\infty$, we obtain the desired result
$$\lim_{t\rightarrow+\infty}\|u(t)-\hat U\|_1=0.$$

\section{Extinction of the overhead: proof of Theorem \protect\ref{th:over}}\label{s:over}

This section deals with the multi-D Burgers equation (flux $f^B$). Remember that $u_+<u_-$ and that the shock is uniformly non-characteristic. For technical reasons, we do not assume anymore a steady shock~: the velocity $v(u_-,u_+)$ may be non-zero.

\bigskip
 
To begin with, we consider an auxiliary datum 
$$a_+(x)=\max\{u_-,a(x)\}.$$
Thanks to the comparison principle, we know that 
$$u(t)=S_ta\le S_ta_+.$$
By assumption, $a_+-u_-=(a-u_-)^+$ is bounded, compactly supported. 

It was remarked in \cite{DSLS} that there exists a uni-triangular matrix $M$ and a constant vector $Z$, both depending upon $u_-$, such that $\tilde v(t,x):=(S_ta_+)(Mx-tZ)-u_-$ is again a solution of the multi-D Burgers equation.  Applying Proposition \ref{p:DSLS} to  $\tilde v$, we infer
$$\|S_ta_+-u_-\|_\infty\le c_d\|a_+-u_-\|_1^\alpha\, t^{-\beta}= c_d\|(a-u_-)^+\|_1^\alpha\, t^{-\beta}$$
where $\beta$ is a positive exponent.
This shows that $\|S_ta_+-u_-\|_\infty\rightarrow0$ as $t\rightarrow+\infty$, and thus $\sup_xu(t)\rightarrow u_-$. The same trick shows that $\inf_xu(t)\rightarrow u_+$.

\bigskip

{\em Velocity estimates.} Because $u_-\le S_ta_+\le\sup a_+$, Proposition \ref{p:supp} tells us that the support of $S_ta_+-u_-$ expands at most at finite velocity. More precisely, if $t,T>0$, then
$${\rm Supp}(S_{t+T}a_+-u_-)\subset{\rm Supp}(S_{T}a_+-u_-)+tC_T,$$
where $C_T$ is the convex hull of
$$\left\{\frac1{s_2-s_1}\,(f(s_2)-f(s_1))\,;\,s_1,s_2\in(u_-,\sup_xS_Ta_+)\right\}\,.$$
From Taylor Formula, the diameter of $C_T$ can be estimated as follow:
\begin{equation}\label{eq:supps}
C_T\subset B(f'(u_-);c_f\|S_Ta_+-u_-\|_\infty),
\end{equation}
where the finite constant $c_f$ depends only upon the second differential of $f$ over $(u_+,u_-+1)$. Roughly speaking,  the overhead travels approximately at the velocity $f'(u_-)$.

Likewise, the shock velocity between $\hat u_-$ and $u_+$ is a perturbation of that from $u_-$ to $u_+$~:
\begin{equation}\label{eq:shvel}
|v(\hat u_-,u_+)-v(u_-,u_+)|\le c_f|\hat u_--u_-|.
\end{equation}
Unless we encounter an ambiguity, we shall denote from now on $v$ and $\hat v$ for both shock velocities, keeping in mind that $v$ is given, but $\hat v$ depends upon our choice of $\hat u_-$.

\bigskip

Let $\eta\in(0,1)$ be a number, small enough that we can apply Proposition \ref{p:autrechoc} to the profile $\hat U$ with the same front as $U$, and with end states $\hat u_-=u_-+\eta$ and $\hat u_+=u_+$. Choose $T>0$ large enough that $\sup_xS_Ta_+<\hat u_-$, which implies $u(T,\cdot)\le\hat u_-$.

According to Theorem \ref{th:struct}, there exists an $s>0$ such that the (compact) support of $u(T)-U(\cdot-Tv)$ is contained in $D_-+Tv+sW$. Notice that the latter also contains $D_-+Tv$, the domain where $\hat U(\cdot-Tv)\equiv\hat u_-$. Thus either $u(T,x)=U(x-Tv)\le\hat U(x-Tv)\le \hat U(x-Tv-sW)$, or $x\in D_-+Tv+sW$, and then $u(T,x)\le\hat u_-=\hat U(x-Tv-sW)$. We infer that $u(T)\le \hat U(\cdot-Tv-sW)$.
By comparison, there follows
\begin{equation}
\label{eq:aplhatU}
u(t+T)\le\hat U(\cdot-Tv-sW-t\hat v)=:\hat U(\cdot-z-t\hat v).
\end{equation}

\bigskip

Coming back to velocities, we compare that of the disturbance $S_ta_+$ of $u_-$ with that of the shock $\hat u_-\mapsto u_+$. We have
$$C_T-\hat v\subset B(f'(u_-)-v;2c_f\eta)=B(F'(u_-);2c_f\eta).$$
Let us use the coordinates $(r,y)\in\R\times H$, introduced in Section \ref{s:NPshocks}. We define a function $g:\R^d\rightarrow\R$ by $g(x)=r-\rho\psi_0(y)$, where $\rho<1$ is the constant expressing the non-characteristicness of the shock front in (\ref{eq:LipPsi}). Because $F'(u_-)\in A(u_-,u_+)^\circ$, we have $r_{F'(u_-)}\ge\psi_0(y_{F'(u_-)})$. And since $F'(u_-)\ne0$, because the shock is non-characteristic, this implies
$g(F'(u_-))>0$. We may choose therefore $\eta>0$ small enough that 
\begin{equation}
\label{eq:chooeta}
\min\left\{g(x)\,;\, x\in B(F'(u_-);2c_f\eta)\right\}=:\alpha>0.
\end{equation}
\begin{lemma}\label{l:assgd}
Let $\eta>0$ be chosen so that Proposition \ref{p:autrechoc} applies, and (\ref{eq:chooeta}) holds true. 

Then for $t>0$ large enough, we have
$${\rm Supp}(S_Ta_+-u_-)+tC_T\subset D_++z+t\hat v.$$
\end{lemma}

\bepr

It suffices to prove that if $K$ is a compact set, then for $t$ large enough,
$$K+t(C_T-\hat v)\subset D_+.$$
It is therefore enough to prove that
$$K+tB(F'(u_-);2c_f\eta)\subset D_+.$$

Recall that $D_+$ is described by $r>\psi(y)$. Applying (\ref{eq:LipPsi}), we see that $D_+$ contains a conical domain $D_\rho$ of equation $g(x)>\psi(0)$. Thus it will be sufficient to have, for $t$ large enough,
$$K+tB(F'(u_-);2c_f\eta)\subset D_\rho.$$
To this end, we evaluate the minimum value of $g$ over the compact set $K+tB(F'(u_-);2c_f\eta)$. Because $\psi_0$ is sub-additive (convex and homogeneous of degree one), we have for every $x$ in this domain
$$g(x)\ge\min_Kg+t\alpha.$$
There remains to choose
$$t>\frac1\alpha\,(\psi(0)-\min_Kg).$$

\enpr

\bigskip

Let $t>0$ be as in Lemma \ref{l:assgd}. Then 
$${\rm Supp}(S_{t+T}a_+-u_-)\subset D_++z+t\hat v.$$
If $u(t+T,x)>u_-$ then necessarily $S_{t+T}a_+(x)>u_-$ and thus $x\in D_++z+t\hat v$. This implies $\hat U(x-z-t\hat v)=u_+$. But since $S_{t+T}a_+\le\hat U(\cdot-z-t\hat v)$, this means $S_{t+T}a_+(x)\le u_+$, a contradiction. Therefore $u(t+T)\le u_-$ everywhere.

\bigskip

The proof that $u(t)\ge u_+$ everywhere for $t$ large enough is similar. Theorem \ref{th:over} is proven.

\appendix

\section{Auxiliary statements}\label{ap}

The following result tells us that the support of bounded solutions propagates at finite velocity.
\begin{prop}\label{p:supp}
Let $b_1,b_2\in L^\infty(\R^d)$ take values in a bounded interval $J$ and be such that ${\rm Supp}(b_2-b_1)$ is compact. Denote $C$ the convex hull of the compact set
$$\left\{\frac1{s_2-s_1}\,(f(s_2)-f(s_1))\,;\,s_1,s_2\in J\right\}\,$$
where the quotient is replaced by $f'(s_1)$ if $s_2=s_1$. Let $K$ be the convex hull of ${\rm Supp}(b_2-b_1)$. Then the support of $S_tb_2-S_tb_1$ is contained in
$K+tC.$
\end{prop}

\bepr

For each direction $\xi\in S_{d-1}$, we consider a half-plane $H=\{x\in\R^d\,;\,x\cdot\xi>\alpha_\xi\}$ separated from $K$. We also denote $c(\xi)=\max\xi\cdot C$. Denoting $u_j(t)=S_tb_j$, we integrate Kru\v{z}kov's inequality
$$\partial_t|u_2-u_1|+{\rm div}({\rm sgn}(u_2-u_1)\,(f(u_2)-f(u_1)))\le0$$
over the domain
$$\left\{(t,x)\in(0,T)\times\R^d\,;\,x\cdot\xi>\alpha_\xi+c(\xi)t\right\}.$$
We obtain
\begin{eqnarray*}
\int_0^T\int_{x\cdot\xi=\alpha_\xi+c(\xi)t}(c(\xi)|u_2-u_1|-{\rm sgn}(u_2-u_1)\,\xi\cdot(f(u_2)-f(u_1)))\,d\Sigma\, dt & & \\
+\sqrt{1+c(\xi)^2\,}\,\int_{x\cdot\xi>\alpha_\xi+c(\xi)T}|u_2-u_1|(T,x)\,dx & \le & 0,
\end{eqnarray*}
where $d\Sigma$ is the $(d-1)$-dimensional Lebesgue measure.
The first integral is non-negative because $u_1,u_2$ take values in $J$, and we conclude that
$$\int_{x\cdot\xi>\alpha_\xi+c(\xi)T}|u_2-u_1|(T,x)\,dx\le0.$$
This implies that the support of $u_2(T)-u_1(T)$ is contained in the half-space 
$$\{x\,;\,x\cdot\xi\le \alpha_\xi+c(\xi)T\}.$$
There remains to take the intersection as $\xi$ runs over the unit sphere.

\enpr

\bigskip

The following dispersion statement is taken from \cite{DSLS} (see Theorem 1.1).
\begin{prop}\label{p:DSLS}
Consider the multi-D Burgers equation (flux $f^B$). Then there exists a universal constant $c_d<\infty$ and positive exponents 
$$\alpha(d)=\frac{2}{d^2+d+2}\,,\qquad\beta(d)=\frac{2d}{d^2+d+2}\,,$$ 
such that the solutions of (\ref{eq:Bur}) with integrable data obey to $\|S_tu_0\|_\infty\le c_d\|u_0\|_1^\alpha\, t^{-\beta}$.
\end{prop}

\end{document}